\documentclass{amsart}
\usepackage{amsfonts}

\usepackage{amscd}

\newtheorem{theorem}{Theorem}[section]

\newtheorem{lemma}[theorem]{Lemma}

\theoremstyle{definition}
\newtheorem{definition}{Definition}[section]

\theoremstyle{remark}

\begin{document}
\title{A new method in Fano geometry}
\author{Ziv Ran}
\address{University of California at Riverside}
\email{ziv@math.ucr.edu}
\author{Herb Clemens}
\address{University of Utah}
\email{clemens@math.utah.edu}
\date{January, 2000}
\maketitle

\begin{abstract}
We give some bounds on the anticanonical degrees of Fano varieties with
Picard number 1 and mild singularities. The proof is based on a study of
positivity properties of sheaves of differential operators on ample line
bundles.
\end{abstract}

\section{Introduction\label{1}}

\subsection{Fano varieties}

The purpose of this paper is to bound the degrees of large classes of Fano
varieties.

\begin{definition}
A \textit{unipolar }$\Bbb{Q}$-\textit{Fano variety} is an $n$-dimensional
complex projective variety $X$ such that

i) $X$ is normal and $\Bbb{Q}$-factorial,

ii) the set of Weil divisors modulo numerical equivalence forms a group
\begin{equation*}
\Bbb{{Z\cdot }\left\{ D_{X}\right\} \cong Z}
\end{equation*}
with $D_{X}$ a Weil divisor,

iii) the $\Bbb{Q}$-Cartier $\left( -K_{X}\right) $ is ample. We write
\begin{equation*}
-K_{X}=i_{X}\left\{ D_{X}\right\} ,
\end{equation*}
for some positive integer $i_{X}$.
\end{definition}

\noindent The positive integer $i_{X}$ is called the \textit{Weil index} of $%
X$. Also we define
\begin{equation*}
t_{X}
\end{equation*}
to be the smallest positive integer such that
\begin{equation*}
t_{X}D_{X}
\end{equation*}
is Cartier. Then, of course,
\begin{equation*}
\left( K_{X}\right) ^{n}=\frac{\left( t_{X}K_{X}\right) ^{n}}{t_{X}^{n}}.
\end{equation*}
By the Appendix to \S 1 of \cite{Re}, the group of Weil divisors is
isomorphic to the set of saturated, torsion-free rank-one sheaves on $X$,
which in turn is the same as the set of rank-one reflexive sheaves. These
sheaves are called \textit{divisorial sheaves}. For example, on a
Cohen-Macaulay, normal variety, the dualizing sheaf is always divisorial.

We let
\begin{equation*}
X^{\prime }
\end{equation*}
denote the smooth points of $X$. Then any divisorial sheaf is equal to the
push-forward of its restriction to any subset of $X$ whose complement has
codimension at least two, in particular from $X^{\prime }$.

\subsection{One-canonical singularities}

\begin{definition}
A normal, Cohen-Macaulay variety $X$ is said to be $1$\textit{-canonical}
if, for any resolution
\begin{equation*}
\epsilon :Y\rightarrow X
\end{equation*}
the differential
\begin{equation*}
d\epsilon :\epsilon ^{*}\left( \Omega _{X}\right) \rightarrow \Omega _{Y}
\end{equation*}
factors through a map
\begin{equation*}
\epsilon ^{*}\left( \Omega _{X}^{\vee \vee }\right) \rightarrow \Omega _{Y},
\end{equation*}
that is, $1$-forms on $X^{\prime }$ lift to holomorhic forms on $Y$.
\end{definition}

Notice that this is the precise analogue with respect to one-forms of the
condition on top forms which defines canonical singularities. Note also that
the condition is automatically satisfied whenever the natural map
\begin{equation*}
\Omega _{X}\rightarrow \Omega _{X}^{\vee \vee }
\end{equation*}
is surjective. It can be shown that any locally finite quotient of a smooth
in codimension $2$ complete intersection is $1$-canonical (see \cite{Ra}).

\subsection{The theorem}

We choose a very ample polarization
\begin{equation*}
H\equiv high\ multiple\ of\ \left\{ D_{X}\right\}
\end{equation*}
and let
\begin{equation}
C\equiv H^{n-1}\subseteq X^{\prime }  \label{1.5}
\end{equation}
be a generic linear section curve of the embedding of $X$ given by $H.$

\begin{theorem}
i) Let $X$ be a unipolar $\Bbb{Q}$-Fano variety with only log-terminal, $1$%
-canonical singularities. Then
\begin{equation*}
\left( -K\right) _{X}^{n}\leq \left( \max .\left\{ \frac{2\left( C\cdot %
\left( -K_{X}\right) \right) }{\mu _{\min .}\left( T_{X}\right) }%
,i_{X}\right\} \right) ^{n}\leq \left( \max .\left\{ \frac{2\left( C\cdot %
\left( -K_{X}\right) \right) }{\mu _{\min .}\left( T_{X}\right) }%
,t_{X}\left( n+1\right) \right\} \right) ^{n}
\end{equation*}
where
\begin{equation*}
\mu _{\min .}\left( T_{X}\right)
\end{equation*}
is the minimum slope of the subquotients in a Harder-Narasimhan filtration
of the tangent bundle $T_{X}$ of $X$ with respect to the polarization $H$
(see \S \ref{4} below),
 and we always have
$$
 \frac{ C\cdot
( -K_{X})  }{\mu _{\min .}( T_{X})}
\leq i_X.
$$

ii) Let $X$ be a unipolar $\Bbb{Q}$-Fano variety. Suppose that $T_{X}$ is
semi-stable. Then
\begin{equation*}
\frac{C\cdot \left( -K_{X}\right) }{\mu _{\min .}\left( T_{X}\right) }=n
\end{equation*}
so that
\begin{equation*}
\left( -K\right) _{X}^{n}\leq \left( \max .\left\{ 2n,t_{X}\left( n+1\right)
\right\} \right) ^{n}.
\end{equation*}
\end{theorem}

In this paper we give a complete proof of this theorem. A slightly more
general result is proved with a slightly different viewpoint in \cite{Ra},
which also contains some ancillary definitions, examples and technical
details, as well as a number of applications. Here our purpose has been to
clarify the main ideas of the proof, and we have thus steered a direct
course to the main result while trying to make the argument comprehensible
to nonexperts. At this point, we suggest that the first-time reader skip
immediately to \S \ref{10} below in order to get an idea of the strategy of
the proof. Suffice it to say here that the proof is based on positivity
properties of sheaves of differential operators and in particular is
completely independent of rational curves and bend-and-break which were used
heavily in earlier approaches to boundedness of Fano varieties (see \cite{Ko}
and references therein).

Both authors would like to thank Paul Burchard, J\'{a}nos Koll\'{a}r, and
Robert Lazarsfeld for helpful communications and especially to thank James
McKernan for several key suggestions and at least one important correction.

\section{Good resolution of $X\label{res}$}

We will need to understand pull-back of divisorial sheaves under ``nice''
resolution of $X$. Let
\begin{equation*}
\epsilon :Y\rightarrow X
\end{equation*}
be a resolution obtained by a succession of ``modifications,'' where a
modification
\begin{equation*}
X_{i+1}\subseteq \Bbb{P}^{n_{i+1}}
\end{equation*}
of
\begin{equation*}
X_{i}\subseteq \Bbb{P}^{n_{i}}
\end{equation*}
is obtained by blowing up $\Bbb{P}^{n_{i}}$ along a smooth center inside the
singular locus of $X_{i}$, then embedding the proper transform $X_{i+1}$ of $%
X_{i}$ into a projective space $\Bbb{P}^{n_{i+1}}$ via a sufficiently high
multiple of the ample divisor
\begin{equation*}
m_{i}\mathcal{O}_{\Bbb{P}^{n_{i}}}\left( 1\right) -F_{i},
\end{equation*}
$F_{i}$ being the exceptional divisor of the blow-up. Repeating this process
as necessary, we arrive at a smooth projective manifold $Y$ such that the
exceptional locus $E=\bigcup\nolimits_{i}E_{i}$ lies over the singular set
of $X$ and is a is a simple-normal-crossing divisor. $Y$ has Neron-Severi
group (modulo numerical equivalence, given by
\begin{equation*}
NS\left( Y\right) =\Bbb{{Z}\cdot \tilde{D}_{X}\oplus \sum\nolimits_{i}{Z}%
\cdot E_{i}}
\end{equation*}
where $\tilde{D}_{X}$ denotes the proper transform of the Weil divisor $%
D_{X} $.

If $M$ is a divisorial sheaf on $X$, then locally at any singular point of $%
X $,
\begin{equation*}
M=I\cdot L
\end{equation*}
with $I$ the ideal sheaf of an effective Weil divisor and $L$ locally free
of rank one. So we define the ``integral-divisorial'' pull-back $\epsilon
_{id.}^{*}M$ of $M$ to be the line bundle on $Y$ given by sections of $%
\epsilon ^{*}L$ whose order along each divisor $B$ is greater than or equal
to
\begin{equation*}
\min \left\{ ord_{B}\left( f\circ \epsilon \right) :f\in I\right\} .
\end{equation*}
along that divisor. If $M$ is Cartier, then pull-back is just pull-back of
line bundles. However, if $M$ is not Cartier, the natural map
\begin{equation*}
m\epsilon _{id.}^{*}M\rightarrow \epsilon ^{*}mM
\end{equation*}
is not necessarily an isomorphism. Later we will need the fact that, for any
morphism
\begin{equation*}
M\rightarrow F
\end{equation*}
for $M$ divisorial and $F$ locally free, there is an induced morphism
\begin{equation*}
\epsilon _{id.}^{*}M\rightarrow \epsilon ^{*}F.
\end{equation*}
We have, for example, for any positive integer $m^{\prime }$,
\begin{equation}
\epsilon _{id.}^{*}m^{\prime }D_{X}=\tilde{D}_{X}+\sum\nolimits_{i}\tilde{a}%
_{i}\left( m^{\prime }\right) E_{i}  \label{1.1}
\end{equation}
where the $\tilde{a}_{i}\left( m^{\prime }\right) $ are integers. On the
other hand, for $M$ divisorial, one has the divisorial pull-back
\begin{equation*}
\epsilon _{div.}^{*}M:=\frac{1}{m}\epsilon ^{*}mM
\end{equation*}
where $mM$ is Cartier. So we have
\begin{equation}
\epsilon _{div.}^{*}D_{X}=\tilde{D}_{X}+\sum\nolimits_{i}\frac{a_{i}}{t_{X}}%
E_{i}  \label{1.2}
\end{equation}
with the $a_{i}$ non-negative integers. Since there is a standard
multiplication map
\begin{equation*}
\left( \epsilon _{id.}^{*}M\right) ^{m}\rightarrow \left( \epsilon
_{div.}^{*}M\right) ^{m}.
\end{equation*}
we have
\begin{equation}
\tilde{a}_{i}\left( m^{\prime }\right) \geq \frac{m^{\prime }a_{i}}{t_{X}}.
\label{1.3}
\end{equation}

By construction, we have that the divisor
\begin{equation*}
m\tilde{D}_{X}+\sum\nolimits_{i}\left( \frac{a_{i}}{t_{X}}-t_{i}\right)
E_{i}.
\end{equation*}
is ample for $0<$ $t_{i}<<1$ and $m>0$. So, by the Kawamata-Viehweg
Vanishing Theorem,
\begin{equation*}
H^{j}\left( -\epsilon ^{*}mt_{X}D_{X}\right) =0
\end{equation*}
for $j\leq dimX$. Also we have
\begin{eqnarray*}
K_{Y} &=&\epsilon _{div.}^{*}K_{X}+\sum\nolimits_{i}b_{i}E_{i} \\
&=&-i_{X}\tilde{D}_{X}+\sum\nolimits_{i}\left( b_{i}-\frac{i_{X}}{t_{X}}%
a_{i}\right) E_{i}.
\end{eqnarray*}
So by duality
\begin{equation}
H^{j}\left( K_{Y}+m\epsilon ^{*}t_{X}D_{X}\right) =H^{j}\left( \left(
mt_{X}-i_{X}\right) \tilde{D}_{X}+\sum\nolimits_{i}\left( b_{i}+\left( m-%
\frac{i_{X}}{t_{X}}\right) a_{i}\right) E_{i}\right) =0  \label{1.4}
\end{equation}
for all $j,m>0$.

\section{Bounding the Weil index\label{2}}

Referring to $\left( \ref{1.4}\right) $%
\begin{eqnarray*}
\chi \left( a\right) &=&h^{0}\left( K_{Y}+\epsilon ^{*}at_{X}D_{X}\right) \\
&=&h^{0}\left( \left( at_{X}-i_{X}\right) \tilde{D}_{X}+\sum\nolimits_{i}%
\left( b_{i}+\left( a-\frac{i_{X}}{t_{X}}\right) a_{i}\right) E_{i}\right)
\end{eqnarray*}
is zero for $1\leq at_{X}\leq \left( i_{X}-1\right) ,$ since, in that case,
the push-forward of a section would be negative along $D_{X}$. To see that $%
\chi \left( a\right) $ is not identically zero, choose large $a$ such that
the pullback of $at_{X}D_{X}$ to $Y$ vanishes to order at least $i_{X}$ on $%
\tilde{D}_{X}$ and order at least $b_{i}$ on each $E_{i}$. Thus
\begin{equation*}
\left( \frac{i_{X}-1}{t_{X}}\right) \leq \deg \chi \left( a\right) \leq n.
\end{equation*}
so that
\begin{eqnarray*}
\frac{i_{X}-1}{t_{X}} &<&n+1, \\
i_{X} &\leq &t_{X}\left( n+1\right) .
\end{eqnarray*}
(This argument in the smooth case is due to Kobayashi-Ochiai.)\smallskip

Referring to \S \ref{1} we therefore have
\begin{equation*}
C\cdot D_{X}=\frac{-C\cdot K_{X}}{i_{X}}\geq \frac{-C\cdot K_{X}}{
t_{X}\left( n+1\right) }.
\end{equation*}

\section{Atiyah class\label{3}}

Given line bundles $L$ and $L^{\prime }$ on $X^{\prime }$, the set of smooth
points of $X$, let
\begin{equation*}
\frak{D}^{n}\left( L,L^{\prime }\right)
\end{equation*}
denote the sheal of holomorphic differential operators of order $\leq n$ on
sections of $L$ with values in sections of $L^{\prime }$. (If $L=L^{\prime }$%
, we simply write $\frak{D}^{n}\left( L\right) $.) The sequence
\begin{equation*}
0\rightarrow \mathcal{O}_{X^{\prime }}\rightarrow \frak{D}^{1}\left(
\mathcal{O}_{X^{\prime }}\right) \rightarrow T_{X^{\prime }}\rightarrow 0
\end{equation*}
splits as a sequence of left $\mathcal{O}_{X^{\prime }}$-modules but not as
a sequence of right $\mathcal{O}_{X^{\prime }}$-modules. In fact, if we
tensor on the right by a line bundle $L^{*}$ to obtain the exact sequence
\begin{equation*}
0\rightarrow L^{*}\rightarrow \frak{D}^{1}\left( L,\mathcal{O}_{X^{\prime
}}\right) \rightarrow T_{X^{\prime }}\otimes L^{*}\rightarrow 0
\end{equation*}
of left $\mathcal{O}_{X^{\prime }}$-modules and then tensoring this last
sequence on the left by $L$, we obtain the exact sequence
\begin{equation*}
0\rightarrow \mathcal{O}\rightarrow \frak{D}_{T_{X^{\prime }}}^{1}\left(
L\right) \rightarrow T_{X^{\prime }}\rightarrow 0.
\end{equation*}
The obstruction to splitting this last sequence (as a sequence of left
modules) is given by taking a meromorphic section $l_{0}$ on $L$ and
splitting the above sequence over the set where $l_{0}\neq 0,\infty $ via
\begin{equation*}
\frac{\partial }{\partial x}\longmapsto \left( l\mapsto l_{0}\frac{\partial
\left( l/l_{0}\right) }{\partial x}\right) .
\end{equation*}
Writing patching data $\left\{ z_{U}\right\} $ for the divisor of $l_{0}$ we
have
\begin{equation*}
z_{U}^{-1}\frac{\partial \left( f\cdot z_{U}\right) }{\partial x}%
-z_{U^{\prime }}^{-1}\frac{\partial \left( f\cdot z_{U^{\prime }}\right) }{%
\partial x}=f\left( \frac{\partial \log z_{U}}{\partial x}-\frac{\partial
\log z_{U}}{\partial x}\right)
\end{equation*}
so that the obstruction to splitting is given by
\begin{equation*}
c_{1}\left( L\right) \in H^{1}\left( \Omega _{X^{\prime }}^{1}\right) .
\end{equation*}

\section{Harder-Narasimhan filtration\label{4}}

Throughout we will deal with vector bundles and sheaves modulo
codimension-two phenomena. Thus ``torsion'' means torsion along a divisor,
``vector bundle'' means locally free through codimension one, etc.

Since the Picard number of $X$ is one and singularities are of codimension
two, there is an unambiguous notion of stability ($H$-stability) of bundles
on $X^{\prime }.$ We will use in an essential way the Harder-Narasimhan
filtration
\begin{equation*}
E_{1}<\ldots <E_{l\left( E\right) -1}<E_{X^{\prime }}
\end{equation*}
of a vector bundle $E$ with
\begin{equation*}
\frac{E_{i}}{E_{i-1}}
\end{equation*}
semi-stable locally free sheaves such that the slopes
\begin{equation*}
\mu _{i}=\frac{c_{1}\left( \frac{E_{i}}{E_{i-1}}\right) \cdot C}{rk\left(
\frac{E_{i}}{E_{i-1}}\right) }
\end{equation*}
form a strictly decreasing sequence whose extremal elements are denoted as
\begin{equation*}
\mu _{\max .}\left( E\right) ,\ \mu _{\min .}\left( E\right)
\end{equation*}
respectively. By results of Mehta-Ramanathan ([MehR]), the above filtration
restricts to a Harder-Narasimhan filtration on a generic curve $C\subseteq
X^{\prime }$ and conversely, that is, a filtration that restricts to a
HN-filtration on generic $C$, is an HN-filtration.

If
\begin{equation*}
E=T_{X^{\prime }},
\end{equation*}
then
\begin{eqnarray*}
\mu _{\max .}\left( T_{X}\right) &=&\mu \left( T_{1}\right) =\frac{
c_{1}\left( T_{1}\right) \cdot C}{t_{1}}\geq \frac{-K_{X}\cdot C}{n} \\
\mu _{\min .}\left( T_{X}\right) &=&\frac{c_{1}\left( T_{X^{\prime
}}/T_{l\left( T_{X^{\prime }}\right) -1}\right) \cdot C}{t_{l\left(
T_{X^{\prime }}\right) }}
\end{eqnarray*}
where $t_{i}=rk\left( T_{i}/T_{i-1}\right) $. Notice that, if $T_{X^{\prime
}}$ is non-negative, then all the $T_{i}$ are integrable since the map
\begin{eqnarray*}
T_{i}\otimes T_{i} &\rightarrow &\frac{T_{X^{\prime }}}{T_{i}} \\
\xi \otimes \eta &\mapsto &\left[ \xi ,\eta \right]
\end{eqnarray*}
is $\mathcal{O}_{X^{\prime }}$-bilinear and
\begin{equation*}
2\mu _{\min .}\left( T_{i}\right) =2\mu \left( \frac{T_{i}}{T_{i+1}}\right)
>\mu _{\max .}\left( \frac{T_{X^{\prime }}}{T_{i}}\right) .
\end{equation*}

Similarly suppose that, for some line bundle $L$, $\frak{D}^{1}\left(
L\right) $ is non-negative. Let
\begin{equation*}
D_{1}\leq \ldots \leq D_{l\left( {D}\right) }=\frak{D}^{1}\left( L\right)
\end{equation*}
be an HN-filtration. Then the map
\begin{eqnarray*}
D_{i}\otimes D_{i} &\rightarrow &\frac{\frak{D}^{1}\left( L\right) }{{D}_{i}}
\\
\xi \otimes \eta &\mapsto &\left[ \xi ,\eta \right]
\end{eqnarray*}
is also zero.

We begin by examining the slopes of the HN-filtration of
\begin{equation*}
\frak{D}^{1}\left( L,\mathcal{O}_{X^{\prime }}\right)
\end{equation*}
the sheaf of first-order differential operators from a line bundle $L$ on $%
X^{\prime }$ to the structure sheaf $\mathcal{O}_{X^{\prime }}$. Again
recalling that all sheaves are taken ``modulo codimension two,'' suppose
that $T^{\prime }$ is a torsion-free subbundle of $T_{X^{\prime }}$. Then
restricting the symbol map
\begin{equation*}
\frak{D}^{1}\left( \mathcal{O}_{X^{\prime }}\right) \rightarrow T_{X^{\prime
}}
\end{equation*}
to the preimage of $T^{\prime }$, we obtain the sequence
\begin{equation*}
0\rightarrow \mathcal{O}_{X^{\prime }}\rightarrow \frak{D}_{T^{\prime
}}^{1}\left( \mathcal{O}_{X^{\prime }}\right) \rightarrow T^{\prime
}\rightarrow 0
\end{equation*}
is exact, so that the sequence
\begin{equation*}
0\rightarrow L^{*}\rightarrow \frak{D}_{T^{\prime }}^{1}\left( L,\mathcal{O}
_{X^{\prime }}\right) \rightarrow T^{\prime }\otimes L^{*}\rightarrow 0
\end{equation*}
obtained by tensoring \textit{on the right} with $L^{*}$ is also exact.

\section{Relative positivity of first order operators\label{5}}

\begin{lemma}
Suppose $T_{X^{\prime }}$ is positive and $L$ is a line bundle with $L\cdot
C\neq 0$. Then
\begin{equation*}
\frak{D}^{1}\left( L\right)
\end{equation*}
is positive, in fact
\begin{equation*}
\mu _{\min .}\left( \frak{D}^{1}\left( L\right) \right) \geq \min \left\{
C\cdot D_{X},\frac{1}{2}\mu _{\min .}\left( T_{X^{\prime }}\right) \right\}
=:b.
\end{equation*}
\end{lemma}

\begin{proof}
Consider the semi-stable quotient
\begin{equation*}
\frak{D}^{1}\left( L\right) \rightarrow \frac{D_{l}}{D_{l-1}}=:F.
\end{equation*}
in an HN-filtration for $\frak{D}_{1}\left( L\right) $. Thus
\begin{equation*}
\mu _{\min .}\left( \frak{D}^{1}\left( L\right) \right) =\mu \left( F\right)
.
\end{equation*}
The composition
\begin{equation}
\mathcal{O}_{X^{\prime }}\rightarrow \frak{D}^{1}\left( L\right) \rightarrow
F  \label{5.1}
\end{equation}
is either zero or injective (through codimension one). If it is zero then we
have a quotient map
\begin{equation*}
T_{X^{\prime }}\rightarrow F
\end{equation*}
so that $\frak{D}_{1}\left( L\right) $ is positive since $T_{X^{\prime }}$
is. In fact, in that case,
\begin{equation*}
\mu _{\min .}\left( \frak{D}^{1}\left( L\right) \right) \geq \mu _{\min
.}\left( T_{X^{\prime }}\right) .
\end{equation*}

If the composition $\left( \ref{5.1}\right) $ is injective, let
\begin{equation*}
M
\end{equation*}
denote the saturation of the image of $\mathcal{O}$ in $F$. If
\begin{equation}
\mathcal{O}\neq M,  \label{5.2}
\end{equation}
then
\begin{equation*}
c_{1}\left( M\right) \geq D_{X}\cdot C
\end{equation*}
so that, by the semi-stability of $F$,
\begin{equation*}
\mu \left( F\right) \geq \frac{1}{2}\mu _{\min .}\left( T_{X^{\prime
}}\right) .
\end{equation*}

If
\begin{equation*}
\mathcal{O}=M,
\end{equation*}
then
\begin{equation*}
\frac{F}{M}
\end{equation*}
is a torsion-free quotient of $F$ and
\begin{eqnarray*}
c_{1}\left( F\right) &=&c_{1}\left( \frac{F}{M}\right) , \\
rkF &=&1+rk\left( \frac{F}{M}\right)
\end{eqnarray*}
so that
\begin{equation*}
\mu \left( F\right) =\frac{c_{1}\left( \frac{F}{M}\right) }{1+rk\left( \frac{%
F}{M}\right) }
\end{equation*}
where
\begin{equation*}
\mu \left( \frac{F}{M}\right) =\frac{c_{1}\left( \frac{F}{M}\right) }{%
rk\left( \frac{F}{M}\right) }\geq \mu _{\min .}\left( T_{X^{\prime }}\right)
.
\end{equation*}
Now, $M\neq F$ by \S \ref{3}, so
\begin{equation*}
\mu _{\min .}\left( \frak{D}^{1}\left( L\right) \right) \geq \frac{rk\left(
\frac{F}{M}\right) }{1+rk\left( \frac{F}{M}\right) }\cdot \mu _{\min
.}\left( T_{X^{\prime }}\right) \geq \frac{1}{2}\mu _{\min .}\left(
T_{X^{\prime }}\right) .
\end{equation*}
\end{proof}

\section{Extending to the case of vector bundles\label{6}}

\begin{lemma}
\label{L6}If $E$ is a positive vector bundle on $X^{\prime }$, then
\begin{equation*}
\mu _{\min .}\left( \frak{{D}^{1}}\left( E,\mathcal{O}\right) \right) \geq
\mu _{\min .}\left( E^{*}\right) +b
\end{equation*}
where
\begin{equation*}
b=\min \left\{ D_{X}\cdot C,\frac{1}{2}\mu _{\min .}\left( T_{X^{\prime
}}\right) \right\} .
\end{equation*}
\end{lemma}

\begin{proof}
Via a HN-filtration for $E$ and the isomorphism
\begin{equation*}
\frak{D}^{1}\left( E^{\prime },\mathcal{O}\right) \rightarrow \frac{\frak{D}%
^{1}\left( E,\mathcal{O}\right) }{\frak{D}^{1}\left( E/E^{\prime },\mathcal{O%
}\right) }
\end{equation*}
reduce to the case $E$ semi-stable. Let
\begin{equation*}
M=\det E.
\end{equation*}
We know from \S \ref{5} that
\begin{equation*}
\mu _{\min .}\left( \frak{{D}^{1}}\left( M,\mathcal{O}\right) \right) \geq
-M\cdot C+b.
\end{equation*}

Next, deform to the normal cone. Namely blow up $\left( C\times \left\{
0\right\} \right) $ in $\left( X\times \Bbb{{A}^{1}}\right) $ and pull $E$
back to the product and let $E^{\prime }$ be the restriction to the normal
bundle $N_{C\backslash X}$ lying inside exceptional divisor. Since the
deformation to the normal cone is trivial on a first-order neighborhood of
the proper transform of $\left( C\times \Bbb{{A}^{1}}\right) ,$ we have
\begin{equation*}
\frak{{D}^{1}}\left( E,\mathcal{O}_{C}\right) =\frak{{D}^{1}}\left(
E^{\prime },\mathcal{O}_{C}\right) =\frak{{D}^{1}}\left( \nu ^{*}E_{C},%
\mathcal{O}_{C}\right)
\end{equation*}
where $\nu :N_{C\backslash X}\rightarrow C$ is the projection given by the
normal bundle. Let $e=rkE$. Taking an unramified $e$-fold cover
\begin{equation*}
\pi :\tilde{C}\rightarrow C,
\end{equation*}
we have
\begin{equation*}
\nu ^{*}E_{C}\times _{C}\tilde{C}=L\otimes F
\end{equation*}
with $L$ the pullback of a line bundle $L_{\tilde{C}}$.on $\tilde{C}$ and $F$
the pull-back of a semi-stable vector bundle $F_{\tilde{C}}$.on $\tilde{C}$.
Also
\begin{eqnarray*}
e\cdot \left( c_{1}L\right) &\equiv &\pi ^{*}\det E \\
c_{1}F &\equiv &0.
\end{eqnarray*}
So by the well-known theorem of Narasimhan-Seshadri, $F_{\tilde{C}}$ is
given by a locally constant sheaf on $\tilde{C}$. Also the product rule
induces an isomorphism
\begin{eqnarray*}
\frak{{D}^{1}}\left( L,\mathcal{O}_{\tilde{C}}\right) &\rightarrow &\frak{{D}%
^{1}}\left( L^{e},L^{e-1}\right) =L^{e-1}\otimes \frak{{D}^{1}}\left( L^{e},%
\mathcal{O}_{\tilde{C}}\right) \\
D &\mapsto &\left( l_{1}\cdot l_{2}\cdot \ldots \cdot l_{e}\mapsto
Dl_{1}\cdot l_{2}\cdot \ldots \cdot l_{e}+l_{1}\cdot Dl_{2}\cdot \ldots
\cdot l_{e}+\ldots \right)
\end{eqnarray*}
so that, since HN-filtrations are preserved under covers, we have by the
rank one case that
\begin{eqnarray}
\mu _{\min .}\left( \frak{{D}^{1}}\left( L,\mathcal{O}_{\tilde{C}}\right)
\right) &=&\mu _{\min .}\left( L^{e-1}\otimes \frak{{D}^{1}}\left( L^{e},%
\mathcal{O}_{\tilde{C}}\right) \right)  \label{6.1} \\
&=&\left( e-1\right) L\cdot \tilde{C}+\mu _{\min .}\left( \pi ^{*}\frak{{D}%
^{1}}\left( \det \left( \nu ^{*}E_{C}\right) ,\mathcal{O}_{N}\right) \right)
\notag \\
&=&\left( e-1\right) L\cdot \tilde{C}+e\cdot \mu _{\min .}\left( \frak{{D}%
^{1}}\left( \det \left( \nu ^{*}E_{C}\right) ,\mathcal{O}_{N}\right) \right)
\notag \\
&=&\left( e-1\right) L\cdot \tilde{C}+e\cdot \mu _{\min .}\left( \frak{{D}%
^{1}}\left( \det E,\mathcal{O}\right) \right)  \notag \\
&\geq &\left( e-1\right) L\cdot \tilde{C}+e\cdot \left( -\det E_{C}+b\right)
\notag \\
&=&-L\cdot \tilde{C}+eb  \notag
\end{eqnarray}
where $N=N_{C\backslash X}$. So, since $E$ is semi-stable
\begin{equation*}
\mu _{\min .}\left( \frak{{D}^{1}}\left( L,\mathcal{O}_{\tilde{C}}\right)
\right) \geq -\frac{\pi ^{*}\det E}{e}+eb=e\left( -\mu \left( E\right)
+b\right) .
\end{equation*}

On the other hand, since $F_{\tilde{C}}$ is locally constant and therefore $%
F $ is too, we have that
\begin{equation*}
\frak{{D}^{1}}\left( F\otimes L,\mathcal{O}_{\pi ^{-1}N}\right) \frak{=}%
F\otimes \frak{{D}^{1}}\left( L,\mathcal{O}_{\pi ^{-1}N}\right)
\end{equation*}
so that
\begin{eqnarray}
\mu _{\min .}\left( \frak{{D}^{1}}\left( E,\mathcal{O}_{C}\right) \right)
&=&\mu _{\min .}\left( \frak{{D}^{1}}\left( \nu ^{*}E_{C},\mathcal{O}%
_{C}\right) \right)  \label{6.2} \\
&=&e^{-1}\cdot \mu _{\min .}\left( \pi ^{*}\frak{{D}^{1}}\left( \det \left(
\nu ^{*}E_{C}\right) ,\mathcal{O}_{\tilde{C}}\right) \right)  \notag \\
&=&e^{-1}\cdot \mu _{\min .}\left( \frak{{D}^{1}}\left( F\otimes L,\mathcal{O%
}_{\tilde{C}}\right) \right)  \notag \\
&=&e^{-1}\cdot \mu _{\min .}\left( F^{*}\otimes \frak{{D}^{1}}\left( L,%
\mathcal{O}_{\tilde{C}}\right) \right)  \notag \\
&=&e^{-1}\cdot \mu _{\min .}\left( \frak{{D}^{1}}\left( L,\mathcal{O}_{%
\tilde{C}}\right) \right) .  \notag
\end{eqnarray}
Taken together $\left( \ref{6.1}\right) $ and $\left( \ref{6.2}\right) $
complete the proof in the case $E$ semi-stable. Therefore we are done.
\end{proof}

\section{Extending to estimates to higher-order operators\label{7}}

In this section we work over $X^{\prime }$, the set of smooth points of $X$.
The extension of the above estimates to higher order operators will be made
using the elementary fact that (restricting to $X^{\prime }$) there is a
natural surjection
\begin{equation}
\frak{D}^{1}\left( P^{i}\left( E\right) ,\mathcal{O}\right) \rightarrow
\frak{D}^{m+1}\left( E,\mathcal{O}\right) .  \label{7.1}
\end{equation}
where $P^{i}\left( E\right) $ is the sheaf of $i$-th order jets of the
vector bundle $E$, that is,
\begin{equation*}
p_{*}q^{*}E
\end{equation*}
where $p$ and $q$ are the two projections of the $i$-th order neighborhood
of the diagonal of $X^{\prime }\times X^{\prime }$. To see this, assigning
to any section its $i$-th jet to get
\begin{equation*}
\frak{D}^{m}\left( E,\mathcal{O}\right) \frak{={Hom}}\left( P^{i}\left(
E\right) ,\mathcal{O}\right)
\end{equation*}
where $\frak{Hom}$ is with reference to the left $\mathcal{O}$-linear
structure. So we have a left $\mathcal{O}$-linear surjection
\begin{equation*}
\frak{{D}^{1}}\left( P^{i}\left( E\right) ,\mathcal{O}\right) =\frak{{D}%
^{1}\left( \mathcal{O}\right) \otimes }Hom\left( P^{i}\left( E\right) ,%
\mathcal{O}\right) =\frak{{D}^{1}}\left( \mathcal{O}\right) \otimes \frak{D}%
^{m}\left( E,\mathcal{O}\right) \rightarrow \frak{D}^{m+1}\left( E,\mathcal{O%
}\right) \frak{.}
\end{equation*}

\begin{lemma}
\label{L7}Suppose that $T_{X^{\prime }}$ is positive so that
\begin{equation*}
b=\min .\left\{ C\cdot D_{X},\frac{1}{2}\mu _{\min .}\left( T_{X}\right)
\right\} >0.
\end{equation*}
If $E$ is a positive vector bundle on $X^{\prime }$,
\begin{equation*}
\mu _{\min .}\left( \frak{D}^{m+1}\left( E,\mathcal{O}\right) \right) \geq
\mu _{\min .}\left( \frak{{D}^{1}}\left( P^{i}\left( E\right) ,\mathcal{O}%
\right) \right) \geq \min .\left\{ 0,\mu _{\min .}\left( E^{*}\right)
+\left( m+1\right) b\right\} .
\end{equation*}
\end{lemma}

\begin{proof}
The first inequality is obtained from $\left( \ref{7.1}\right) $. For the
second, we proceed by induction on $m$. The case $m=0$ comes from Lemma \ref
{L6}. Suppose now that, for the quotient
\begin{equation*}
Q
\end{equation*}
of minimal slope in a HN-filtration of
\begin{equation*}
\frak{D}^{m}\left( E,\mathcal{O}\right) \frak{,}
\end{equation*}
we know that
\begin{equation*}
\mu \left( Q\right) \geq \min .\left\{ 0,\left( \mu _{\min .}\left(
E^{*}\right) +mb\right) \right\} .
\end{equation*}
Now if
\begin{equation*}
\mu \left( Q\right) <0,
\end{equation*}
we have that the semi-stable bundle $Q^{*}$ is positive, so by Lemma \ref{L7}%
,
\begin{equation*}
\mu _{\min .}\left( \frak{D}^{1}\left( Q^{*},\mathcal{O}\right) \right) \geq
\mu \left( Q\right) +b\geq \min .\left\{ b,\left( \mu _{\min .}\left(
E^{*}\right) +\left( m+1\right) b\right) \right\} .
\end{equation*}
On the other hand, if
\begin{equation*}
\mu \left( Q\right) \geq 0,
\end{equation*}
then the sequence
\begin{equation*}
0\rightarrow Q\rightarrow \frak{D}^{1}\left( Q^{*},\mathcal{O}\right)
\rightarrow Q\otimes T_{X^{\prime }}\rightarrow 0,
\end{equation*}
the positivity of the tangent bundle, and the nice behavior of
semi-stability under tensor product gives that
\begin{equation*}
\mu _{\min .}\left( \frak{D}^{1}\left( Q^{*},\mathcal{O}\right) \right) \geq
0.
\end{equation*}

Now, dualizing the exact sequence
\begin{equation*}
0\rightarrow S\rightarrow \frak{D}^{m}\left( E,\mathcal{O}\right)
\rightarrow Q\rightarrow 0,
\end{equation*}
gives
\begin{equation*}
0\rightarrow Q^{*}\rightarrow P^{m}\left( E\right) \rightarrow
S^{*}\rightarrow 0
\end{equation*}
from which we obtain the isomorphism
\begin{equation*}
\frak{D}^{1}\left( Q^{*},\mathcal{O}\right) \rightarrow \frac{\frak{D}%
^{1}\left( P^{m}\left( E\right) ,\mathcal{O}\right) }{\frak{D}^{1}\left(
S^{*},\mathcal{O}\right) }\frak{.}
\end{equation*}
Thus
\begin{equation*}
\mu _{\min .}\left( \frak{D}^{1}\left( P^{m}\left( E\right) ,\mathcal{O}%
\right) \right) \geq \mu _{\min .}\left( \frak{D}^{1}\left( Q^{*},\mathcal{O}%
\right) \right)
\end{equation*}
which completes the proof.
\end{proof}

\section{Asymptotic semi-positivity}

Specializing Lemma \ref{L7} to the case of line bundles $L$, we have
\begin{equation*}
\mu _{\min .}\left( \frak{D}^{1}\left( P^{i}\left( L\right) ,\mathcal{O}%
\right) \right) \geq \min .\left\{ 0,-L\cdot C+\left( m+1\right) b\right\}
\end{equation*}
from which we immediately conclude:

\begin{lemma}
If $T_{X}$ is positive and $L$ is any line bundle on $X$ and
\begin{equation*}
\left( m+1\right) \geq \frac{L\cdot C}{b},
\end{equation*}
then $\frak{{D}^{1}}\left( P^{m}\left( L\right) ,\mathcal{O}\right) $ and so
also $\frak{D}^{m+1}\left( L,\mathcal{O}\right) $ are semi-positive.
\end{lemma}

So if
\begin{equation*}
\frac{C\cdot \left( -K_{X}\right) }{b}k\leq \left( m+1\right)
\end{equation*}
we have that
\begin{equation}
\frak{D}^{1}\left( P^{m}\left( -kK_{X}\right) ,\mathcal{O}\right) ,\frak{D}%
^{m+1}\left( \left( -kK_{X}\right) ,\mathcal{O}\right)  \label{8.1}
\end{equation}
are both semipositive. Recall that
\begin{equation*}
b=\min .\left\{ C\cdot D_{X},\frac{1}{2}\mu _{\min .}\left( T_{X}\right)
\right\}
\end{equation*}
so that
\begin{equation*}
\frac{C\cdot \left( -kK_{X}\right) }{b}=\max .\left\{ \frac{2\left( C\cdot
\left( -kK_{X}\right) \right) }{\mu _{\min .}\left( T_{X}\right) }%
,ki_{X}\right\} .
\end{equation*}
So if
\begin{equation*}
\alpha \geq \max .\left\{ \frac{2\left( C\cdot \left( -K_{X}\right) \right)
}{\mu _{\min .}\left( T_{X}\right) },i_{X}\right\} ,
\end{equation*}
then whenever
\begin{equation*}
\alpha k
\end{equation*}
is an integer we have that
\begin{equation*}
\frak{D}^{\alpha k}\left( \left( -kK_{X}\right) ,\mathcal{O}\right)
\end{equation*}
is semipositive.

\begin{lemma}
If
\begin{equation*}
\alpha \geq \max .\left\{ \frac{2\left( C\cdot \left( -K_{X}\right) \right)
}{\mu _{\min .}\left( T_{X}\right) },i_{X}\right\} ,
\end{equation*}
then whenever
\begin{equation*}
\alpha k
\end{equation*}
is a sufficiently divisible integer,
\begin{equation*}
\frak{D}^{\alpha k}\left( -kK_{X},\mathcal{O}_{C}\right)
\end{equation*}
is semi-positive (for sufficiently general $C$).
\end{lemma}

\section{Positivity of the tangent bundle}

>From this point on we restrict the (normal) singularities we allow on $X$.
The necessity of considering only log terminal $X$ derives from the
following:

\begin{lemma}
If $X$ is a log-terminal, $1$-canonical unipolar $\Bbb{Q}$-Fano variety,
then
\begin{equation*}
T_{X}
\end{equation*}
is positive, that is, it has no quotient $Q$ which is locally free in
codimension one and has non-positive first Chern class
\begin{equation*}
c_{1}\left( Q\right) \in \left( \Bbb{Z-N}\right) D_{X}.
\end{equation*}
\end{lemma}

\begin{proof}
Let
\begin{equation*}
\epsilon :Y\rightarrow X
\end{equation*}
be as in \ref{res}. Then using $\left( \ref{1.1}\right) $-$\left( \ref{1.2}%
\right) $%
\begin{equation*}
\epsilon _{id.}^{*}m^{\prime }D_{X}=\tilde{D}_{X}+\sum\nolimits_{i}\tilde{a}%
_{i}\left( m^{\prime }\right) E_{i}
\end{equation*}
and
\begin{equation*}
\epsilon _{div.}^{*}D_{X}=\tilde{D}_{X}+\sum\nolimits_{i}\frac{a_{i}}{t_{X}}%
E_{i}
\end{equation*}
There is an ample $\Bbb{Q}$-divisor
\begin{eqnarray*}
A &:&=\epsilon _{div.}^{*}D_{X}-\sum\nolimits_{i}t_{i}E_{i} \\
&=&\tilde{D}_{X}+\sum\nolimits_{i}\left( \frac{a_{i}}{t_{X}}-t_{i}\right)
E_{i}
\end{eqnarray*}
with $0<t_{i}<<1$. Suppose $Q$ is a torsion-free quotient of
\begin{equation*}
T_{X}
\end{equation*}
with
\begin{equation*}
c_{1}\left( Q\right) =-m^{\prime }D_{X},\ m^{\prime }\geq 0.
\end{equation*}
Then
\begin{equation*}
rk\ Q=n^{\prime }<n
\end{equation*}
since $\bigwedge\nolimits^{n}T_{X}$ is positive. Define
\begin{equation*}
M:=\left( \bigwedge\nolimits^{n^{\prime }}Q^{\vee }\right) ,
\end{equation*}
and consider the natural map
\begin{equation*}
M\rightarrow \bigwedge\nolimits^{n^{\prime }}\left( \Omega _{X}^{\vee \vee
}\right)
\end{equation*}
Then $M^{\vee \vee }$is a divisorial sheaf on $X$, so, using $\left( \ref
{1.3}\right) $%
\begin{eqnarray*}
c_{1}\left( \epsilon _{id.}^{*}M^{\vee \vee }\right) &=&m^{\prime }\tilde{D}%
_{X}+\sum\nolimits_{i}\tilde{a}_{i}\left( m^{\prime }\right) E_{i}. \\
m^{\prime }A &=&m^{\prime }\tilde{D}_{X}+\sum\nolimits_{i}\left( \frac{%
m^{\prime }a_{i}}{t_{X}}-t_{i}\right) E_{i} \\
c_{1}\left( \epsilon _{id.}^{*}M^{\vee \vee }\right) -m^{\prime }A
&=&\sum\nolimits_{i}\left( \tilde{a}_{i}\left( m^{\prime }\right) -\frac{%
m^{\prime }a_{i}}{t_{X}}+t_{i}\right) E_{i} \\
&=&B+\sum\nolimits_{i}s_{i}E_{i}
\end{eqnarray*}
with $B$ integral, effective and $0\leq s_{i}<1$. Thus one can write
\begin{equation}
c_{1}\left( \epsilon _{id.}^{*}M^{\vee \vee }\right) =m^{\prime
}A+B+\sum\nolimits_{i}s_{i}E_{i}  \label{9.1}
\end{equation}
while $A$ is the ample divisor given above, and $0\leq s_{i}<1$.

On the other hand, the map
\begin{equation*}
Q^{\vee }\rightarrow \Omega _{X}^{\vee \vee }
\end{equation*}
and so, by the $1$-canonical condition, induces a map
\begin{equation*}
\epsilon ^{*}Q^{\vee }\rightarrow \epsilon ^{*}\left( \Omega _{X}^{\vee \vee
}\right) \rightarrow \Omega _{Y}
\end{equation*}
and so one gets a non-trivial map
\begin{equation*}
\epsilon ^{*}\left( M^{\vee \vee }\right) \rightarrow \Omega _{Y}^{n^{\prime
}}
\end{equation*}
and therefore a natural map
\begin{equation*}
N:=\epsilon _{id.}^{*}\left( M^{\vee \vee }\right) \rightarrow \Omega
_{Y}^{n^{\prime }}.
\end{equation*}
Thus one must have $H^{0}\left( \Omega _{Y}^{n^{\prime }}\left( -N\right)
\right) \neq 0.$

To contradict the existence of $Q$, we show that
\begin{equation}
H^{0}\left( \Omega _{Y}^{n^{\prime }}\left( -N+B\right) \right) =0.
\label{9.2}
\end{equation}

Case One: $m^{\prime }=0$.

If $M^{\vee \vee }=\mathcal{O}_{X}$, then $N=\mathcal{O}_{Y}$, and one must
show $H^{n^{\prime }}\left( \mathcal{O}_{Y}\right) =0$. By log-terminal
Kodaira Vanishing (Theorem 1-2-5 of \cite{KMM}), $H^{j}\left( \mathcal{O}%
_{X}\right) =0$ for all $j>0$, and by the rationality of log-terminal
singularities (Theorem 1-3-6 of \cite{KMM}), the higher direct image sheaves
$R^{j}\epsilon _{*}\left( \mathcal{O}_{Y}\right) =0$ for all $j>0.$ So by
the Leray spectral sequence, $H^{n^{\prime }}\left( \mathcal{O}_{Y}\right)
=0 $. If $M^{\vee \vee }$ and hence $N$ is only numerically trivial, then
the fact that $H^{1}\left( \mathcal{O}_{Y}\right) =0$ (by the argument just
above) implies that $N=\mathcal{O}_{Y}$.

Case Two: $m^{\prime }>0$.

In this case, one shows $\left( \ref{9.2}\right) $ by using the
branched-covering technique employed in the proof of the Kawamata-Viehweg
Vanishing Theorem. Namely use Theorem 1-1-1 of \cite{KMM} to construct a
smooth finite Galois cover
\begin{equation*}
\tau :Z\rightarrow Y
\end{equation*}
for which the ample $\Bbb{Q}$-divisor
\begin{equation*}
\tau ^{*}A
\end{equation*}
is actually integral (see $\left( \ref{9.1}\right) $). Then
\begin{equation*}
\Omega _{Y}^{n^{\prime }}\left( -N\right) =\Omega _{Y}^{n^{\prime }}\left(
-A-\sum\nolimits_{i}s_{i}E_{i}\right)
\end{equation*}
is a subsheaf of
\begin{equation*}
\tau _{*}\left( \tau ^{*}\Omega _{Y}^{n^{\prime }}\left( -A\right) \right) .
\end{equation*}
Since we have an injection
\begin{equation*}
\tau ^{*}\Omega _{Y}^{n^{\prime }}\rightarrow \Omega _{Z}^{n^{\prime }},
\end{equation*}
and since $\tau ^{*}A$ is ample,
\begin{equation*}
H^{0}\left( \Omega _{Z}^{n^{\prime }}\left( -\tau ^{*}A\right) \right) =0
\end{equation*}
by the Nakano Vanishing Theorem, and so
\begin{equation*}
H^{0}\left( \tau _{*}\left( \tau ^{*}\Omega _{Y}^{n^{\prime }}\left(
-A\right) \right) \right) =0,
\end{equation*}
which completes the proof.
\end{proof}

\section{The strategy/completion of the proof of the Theorem\label{10}}

\subsection{The assumption}

Let
\begin{equation*}
\mu _{X}=\max .\left\{ C\cdot D_{X},\frac{2C\cdot \left( -K_{X}\right) }{\mu
_{\min .}\left( T_{X}\right) }\right\}
\end{equation*}
and suppose
\begin{equation*}
\left( -K_{X}\right) ^{n}>\mu _{X}^{n}.
\end{equation*}
Choose rational constants $\alpha ,\beta $ with
\begin{equation*}
\left( -K_{X}\right) ^{n}>\beta ^{n}>\alpha ^{n}>\mu _{X}^{n}.
\end{equation*}

\subsection{Asymptotic lower bound on sections\label{10.2}}

For sufficiently divisible $k\in \Bbb{N}$, the Hilbert polynomial
\begin{equation*}
\chi \left( -kK_{X}\right) =\frac{\left( -K_{X}\right) ^{n}}{n!}k^{n}+\left(
lower\ powers\ of\ k\right)
\end{equation*}
gives an asymptotic lower bound on
\begin{equation*}
h^{0}\left( -kK_{X}\right) .
\end{equation*}
It is $\leq \binom{m+n-1}{n}$ conditions that a section $s$ of $\left(
-kK_{X}\right) $ have a zero of order $m$ at a determined point $x_{0}\in X$.

\subsection{Semipositivity\label{10.3}}

If $k\alpha $ is an integer and $k$ is sufficiently divisible, the sheaf of
differential operators
\begin{equation*}
\frak{D}^{\alpha k}\left( -kK_{X},\mathcal{O}_{X}\right)
\end{equation*}
was shown in \S{9} to be a semipositive bundle, that is, for any
sufficiently general complete curve-section
\begin{equation*}
C\subseteq X^{\prime }
\end{equation*}
the vector bundle
\begin{equation*}
\frak{D}^{\alpha k}\left( -kK_{X},\mathcal{O}_{C}\right) ={Hom}\left(
\mathcal{O}_{X},\mathcal{O}_{C}\right) \otimes _{\mathcal{O}_{X}} {D}%
^{\alpha k}\left( -kK_{X},\mathcal{O}_{X}\right)
\end{equation*}
has no quotients of negative degree.

\subsection{The contradiction}

For $k>>0$, suppose that, for the $n$-th degree polynomial
\begin{equation*}
\binom{m+n-1}{n}+1=\frac{1}{n!}m^{n}+\ldots ,
\end{equation*}
we let
\begin{equation*}
m_{k}=\alpha k+1
\end{equation*}
so that
\begin{equation*}
\beta k>m_{k}>m_{k}-1=\alpha k.
\end{equation*}
Then
\begin{equation*}
h^{0}\left( -kK_{X}\right) >\frac{\beta ^{n}}{n!}k^{n}\geq \binom{m_{k}+n-1}{%
n}+1.
\end{equation*}
So by \ref{10.2} we have a non-trivial
\begin{equation*}
s\in h^{0}\left( -kK_{X}\right)
\end{equation*}
with a zero of order at least $\alpha k+1$ at a given point $x\in C$, our
general curve. But the mapping
\begin{eqnarray*}
\frak{D}^{\alpha k}\left( -kK_{X},\mathcal{O}_{C}\right) &\longrightarrow &%
\mathcal{O}_{C} \\
D &\longmapsto &D\left( s\right)
\end{eqnarray*}
cannot be trivial because no function in a neighborhood of $x$ in $X$ is
annihilated by all differential operator of degree $\leq \alpha k$. On the
other hand, this last map factors through
\begin{equation*}
\mathcal{O}_{C}\left( -x\right) ,
\end{equation*}
a negative line bundle, contradicting \ref{10.3}. \vfill\eject


\begin{thebibliography}{MehR}
\bibitem[KMM]{KMM}  Kawamata, K., Matsuda, K., Matsuki, K. ``Introduction to
the minimal model program.'' \textit{Adv. Studies Pure Math.} \textbf{10}%
(1987), 283-360.

\bibitem[Ko]{Ko}  Koll\'{a}r, J. ``Rational curves on algebraic varieties.''
Springer 1997.

\bibitem[MehR]{Mehr}  Mehta, V., Ramanathan, A.: ``Semi-stable sheaves on
projective varieties and their restriction to curves''. Math. Ann. \textbf{\
258} (1982), 213-224.

\bibitem[Ra]{Ra}  Ran, Z. ``On semi-positivity of sheaves of differential
operators and the degree of a unipolar $\Bbb{Q}$-Fano variety.'' Preprint,
University of California at Riverside (1999)(eprint math.AG/9811022).

\bibitem[Re]{Re}  Reid, M. ``Canonical $3$-folds.'' \textit{Journ\`{e}es
d'Angers.} North-Holland(1980).
\end{thebibliography}
\end{document}